\documentclass[final,12pt]{elsarticle}
\usepackage{ragged2e}
\justifying
\usepackage{bm}
\usepackage{geometry}
\usepackage{amsmath}
\usepackage{mathrsfs}
\usepackage{amsthm}
\usepackage{amssymb}
\newtheorem{theorem}{Theorem}[section]

\newtheorem{lemma}{Lemma}[section]

\newtheorem{conjecture}{Conjecture}[section]
\newtheorem{proposition}{Proposition}[section]

\numberwithin{equation}{section}
\newenvironment{proof1}{\noindent{\textbf{Proof.}}\ }{\hfill $\square$\par}

\usepackage{indentfirst}
\begin{document}
\begin{frontmatter}
  \title{Some signed graphs whose eigenvalues are main\,\tnoteref{titlenote}}
  \tnotetext[titlenote]{This work was supported by the National Nature Science Foundation
   of China (Nos. 11871040).}
  
  \author{Zhenan Shao}
  \ead{shaozhenan@gmail.com}
 
  \author{Xiying Yuan\corref{correspondingauthor}}
  \cortext[correspondingauthor]{Corresponding author. Email: {\tt xiyingyuan@shu.edu.cn} (Xiying Yuan).}
  
  \address{Department of Mathematics, Shanghai University, Shanghai 200444, P.R. China}
      
  \begin{abstract}
Let $G$ be a graph. For a subset $X$ of $V(G)$, the switching $\sigma$ of $G$ is the signed graph $G^{\sigma}$ obtained from $G$ by reversing the signs of all edges between $X$ and $V(G)\setminus X$. Let $A(G^{\sigma})$ be the adjacency matrix of $G^{\sigma}$. An eigenvalue of $A(G^{\sigma})$ is called a main eigenvalue if it has an eigenvector the sum of whose entries is not equal to zero. Let $S_{n,k}$ be the graph obtained from the complete graph $K_{n-r}$ by attaching $r$ pendent edges at some vertex of $K_{n-r}$. In this paper we prove that there exists a switching $\sigma$ such that all eigenvalues of $G^{\sigma}$ are main when $G$ is a complete multipartite graph, or $G$ is a harmonic tree, or $G$ is $S_{n,k}$. These results partly confirm a conjecture of 
Akbari et al.
 \end{abstract}
  
 \end{frontmatter}
\section{Introduction}
\indent Let $G$ be a graph with vertex set $V(G)$ and edge set $E(G)$. A signed graph $\Gamma=(G,\tau)$ consists of the underlying graph $G$ and a function $\tau$ defined on the edge set $E(G)$, $\tau:E(G)\to\{+1,-1\}$. The graph $G$ can be interpreted as a signed graph in which all edges are positive. For a signed graph $\Gamma$, the \emph{adjacency matrix} of $\Gamma$ is $A(\Gamma)= (a_{ij})$ with $a_{ij}=1$ (resp. $-1$) if $v_{i}$ and $v_{j}$ are adjacent with a positive (resp. negative) edge and $a_{ij}=0$ otherwise. The eigenvalues of $A(\Gamma)$ are called the $\emph{eigenvalues}$ of $\Gamma$ in this paper. We use $\lambda^{[m]}$ to denote an eigenvalue $\lambda$ of $\Gamma$ with multiplicity $m$, and $\lambda^{[0]}$ means $\lambda$ is not an eigenvalue of $\Gamma$. Obviously, $A(\Gamma)$ is a real symmetric matrix and hence all eigenvalues are real. An eigenvalue of $A(\Gamma)$ is called a \emph{main eigenvalue}, if it has an eigenvector the sum of whose entries is not equal to zero. Such a vector is called a $\emph{main eigenvector}$. Let $N(u)$ be the neighborhood of a vertex $u\in V(G)$. The \emph{closed neighbors set} of $u$ is $N[u]=N(u)\cup\{u\}$. Let $\bm{e_{i}}$ be the vector with all-0 entries except for $i$-th entry. Let $\bm{j_{n}}$ be the vector with all-1 entries, and we usually use $\bm j$ instead of $\bm{j_{n}}$. \\
\indent For a subset $X$ of $V(G)$, the \emph{switching} $\sigma$ of $G$ about $X$ is the signed graph $G^{\sigma}$ obtained from $G$ by reversing the signs of all the edges between $X$ and $V(G)\backslash X$. Let $P$=diag$(p_{1},p_{2},\cdots,$\\
$p_{n})$ be a diagonal matrix with $p_{i}=-1$ for the switched vertex $v_{i}$ and 1 otherwise in the switching $\sigma$, then $A(G^{\sigma})=PA(G)P^{-1}$. Hence, $A(G)$ and $A(G^{\sigma})$ have the same spectrum. Moreover, if $\bm u$ is an eigenvector of $A(G)$ associated with the eigenvalue $\lambda$, then $P\bm u$ is an eigenvector of $A(G^{\sigma})$ associated with the eigenvalue $\lambda$.\\
\indent \rm Cvetkovi\'{c} introduced the definition of the main eigenvalue of a graph and proposed the long-standing problem which graphs have exactly $k$ main eigenvalues \cite{dcve}. Over the past few decades,  this problem has been deeply studied. Unicyclic graphs, bicyclic graphs and tricyclic graphs with exactly two main eigenvalues were completely characterized in \cite{hou,hou1}. Huang et al. \cite{huang} gave a construction of graphs with exactly $k$ main eigenvalues and Du et al. \cite{du} constructed infinitely many graphs of order $n$ with $n-1$ main eigenvalues. In recent years, the research on the spectrum of signed graphs has drawn more and more attention \cite{akba,bel,kol,zal}. A natural problem is which signed graphs have exactly $k$ main eigenvalues. Some conclusions about the main eigenvalues of graphs have been generalized to signed graphs. For example, a graph with exactly one main eigenvalue if and only if it is regular \cite{dcve}, and a signed graph with exactly one main eigenvalue if and only if it is net-regular \cite{stan}.  \\
\indent In this paper we first introduce two specific signed graph whose eigenvalues are all main. Denote by $K_{1,n-1}$ the star of order $n$, and by $D_{a,b}$ the double star of order $a+b+2$. For an interger $a\geq2$, the harmonic tree $T_{a}$ is obtained from the tree with one vertex $v$ of degree $a^{2}-a+1$, every neighbour of $v$ with degree $a$ and all remaining vertices with degree one \cite{gru}. For trees, there is no such tree with exactly one main eigenvalue except $K_{1,1}$. The trees with exactly two main eigenvalues are stars, double stars and harmonic trees \cite{hou2}. Let $K_{1,n-1}^{\sigma^{*}}$ be the signed star obtained from $K_{1,n-1}$ $(n\geq3)$ by switching one pendent vertex and the vertex with degree $n-1$. It was proved that all eigenvalues of $K_{1,n-1}^{\sigma^{*}}$ are main \cite{akb}. There exists a switching $\sigma$ such that all eigenvalues of $D_{a,b}^{\sigma}$ are main \cite{akb}. Let $T_{a}^{\sigma^{*}}$ be obtained from $T_{a}$ by switching one pendent vertex and the vertex $v$. We will show that all eigenvalues of $T_{a}^{\sigma^{*}}$ $(a\geq 2)$ are main in Section 2. We use the notation $K_{l_{1}*t_{1},\cdots,l_{s}*t_{s}}$ to denote the complete multipartite graph, which contains $l_{i}$ parts with partition number $t_{i}$, where $t_{1}>t_{2}>\cdots>t_{s}$. For convenience, we use $K_{t_{1}, \cdots,t_{s}}$ instead of $K_{l_{1}*t_{1},\cdots,l_{s}*t_{s}}$ if $l_{i}=1$ for $1\leq i\leq s$. When $t_{s}\geq 2$ and $s\geq 2$ hold, we denote by $K_{t_{1}, \cdots,t_{s}}^{\sigma^{*}}$  the signed complete multipartite graph obtained from $K_{t_{1}, \cdots,t_{s}}$ by switching one vertex in each part. We will also show that all eigenvalues of $K_{t_{1}, \cdots,t_{s}}^{\sigma^{*}}$  are main in Section 2. \\
\indent Graphs $G^{\sigma}$ and $G$ have the same spectrum, while the number of main eigenvalues may be different. Very recently, Akbari et al \cite{akb} proposed the following conjecture, which says there always exists a switching $\sigma$ such that all eigenvalues of $G^{\sigma}$ are main unless $G=K_{2}$ or $K_{4}\backslash\{e\}$. 
\begin{conjecture}\label{conject}
Let $G\neq K_{2},K_{4}\backslash\{e\}$. Then there exists a switching $\sigma$ such that all eigenvalues of $G^{\sigma}$ are main.
\end{conjecture}
\indent Akbari et al. \cite{akb} showed Conjecture 1.1 holds for Cayley graphs, distance-regular graphs, vertex-transitive graphs and edge-transitive regular graphs. In addition, they proved the conjecture holds for double stars, paths and the complete bipartite graphs. Write $S_{n,r}$ for the graph obtained from the complete graph $K_{n-r}$ by attaching $r$ pendent edges at some vertex of $K_{n-r}$. We prove that Conjecture \ref{conject} holds for $S_{n,r}$ in Section 4 and Conjecture \ref{conject} holds for the complete multipartite graph in Section 5. 
\section{All eigenvalues of $T_{a}^{\sigma^{*}}$ and $K_{t_{1}, \cdots,t_{s}}^{\sigma^{*}}$ are main}
\indent Let $\lambda$ be an eigenvalue of $A(\Gamma)$ with multiplicity $m$. If $A(\Gamma)$ has a main eigenvector associated with $\lambda$, then we may construct $m$ main eigenvectors for $\lambda$. In fact, we suppose $\bm{x_{1}},\bm{x_{2}},\cdots,\bm{x_{m}}$ are $m$ independent eigenvectors for $\lambda$ and $\bm y$ is a main eigenvector for $\lambda$. If some $\bm{x_{i}}$ is non-main, we may replace $\bm{x_{i}}$ with $\bm{y_{i}}=\bm y +\bm{x_{i}}$ to obtain a main eigenvector. Hence, we only need to verify that all distinct eigenvalues of $\Gamma$ are main. We always use $\bm x^{t}\cdot\bm y$ to denote the inner product of vectors $\bm x$ and $\bm y$ in this paper. 
\begin{lemma}(\cite[Lemma 3]{hou2})\label{5}
The eigenvalues of $T_{a}$ are $a$, $-a$, $\sqrt{a-1}^{[a^{2}-a]}$, $(-\sqrt{a-1})^{[a^{2}-a]}$ and $0^{[(a-2)(a^{2}-a+1)+1]}$.
\end{lemma}
\begin{proposition}
All eigenvalues of $T_{a}^{\sigma^{*}}$ $(a\geq 2)$ are main.
\end{proposition}
\begin{proof1}
Write $s=(a^{2}-a+1)(a-1)$, $b_{k}=k(a-1)$ and $t=(a^{2}-a+1)a+1$. We label the pendent vertices of $T_{a}$ as $v_{1},v_{2},\cdots, v_{s}$, the vertices adjacent to $v_{b_{k-1}+1},v_{b_{k-1}+2},\cdots,v_{b_{k}}$ as $v_{s+k}$ and the vertices with degree $a^{2}-a+1$ as $v_{t}$. Without loss of generality, we suppose $T_{a}^{*}$ is obtained from $T_{a}$ by switching the vertices $v_{1}$ and $v_{t}$.\\
\indent By Lemma \ref{5}, we know there are five distinct eigenvalues of $T_{a}$ if $a\geq 2$. It may be checked that 
\begin{align*}
&\bm x=\big[\bm{j_{s}}, a\bm{j_{t-s-1}}, a^{2}-a+1\big]^{t},\\
&\bm y=\big[\bm{j_{s}}, -a\bm{j_{t-s-1}}, a^{2}-a+1\big]^{t},\\
&\bm z=\big[\bm{j_{a-1}},-\bm{j_{a-1}},0\bm{j_{s-2a+2}},\sqrt{a-1},-\sqrt{a-1},0\bm{j_{t-s-2}}\big]^{t},\\
&\bm w=\big[\bm{j_{a-1}},-\bm{j_{a-1}},0\bm{j_{s-2a+2}},-\sqrt{a-1},\sqrt{a-1},0\bm{j_{t-s-2}}\big]^{t}
\end{align*}
are eigenvectors associated with the eigenvalue $a$, $-a$, $\sqrt{a-1}$, $-\sqrt{a-1}$, respectively. Set $P$=diag$(-1,1,\cdots,1,-1)$. We have $P\bm x,P\bm y,P\bm z,P\bm w$ are eigenvectors associated with the eigenvalue $a$, $-a$, $\sqrt{a-1}$, $-\sqrt{a-1}$, respectively, and 
\begin{align*}
&P\bm x=\big[-1,\bm{j_{s-1}}, a\bm{j_{t-s-1}}, -a^{2}+a-1\big]^{t},\\
&P\bm y=\big[-1,\bm{j_{s-1}}, -a\bm{j_{t-s-1}}, -a^{2}+a-1\big]^{t},\\
&P\bm z=\big[-1,\bm{j_{a-2}},-\bm{j_{a-1}},0\bm{j_{s-2a+2}},\sqrt{a-1},-\sqrt{a-1},0\bm{j_{t-s-2}}\big]^{t},\\
&P\bm w=\big[-1,\bm{j_{a-2}},-\bm{j_{a-1}},0\bm{j_{s-2a+2}},-\sqrt{a-1},\sqrt{a-1},0\bm{j_{t-s-2}}\big]^{t}.
\end{align*}
Then 
\begin{align*}
&(P\bm x)^{t}\cdot\bm j=2(a^{3}-2a^{2}+2a-2)\neq 0,\\
&(P\bm y)^{t}\cdot\bm j=-2(a^{2}-a+2)\neq 0,\\
&(P\bm z)^{t}\cdot\bm j=-2,\\
&(P\bm w)^{t}\cdot\bm j=-2.
\end{align*}
It is not difficult to know $\bm{e_{1}}+\bm{e_{2}}$ is a main eigenvector associated with the eigenvalue 0 of $T_{a}^{\sigma^{*}}$ when $a\geq 3$. If $a=2$, then $\big[1,-1,-1,0,0,0,-1\big]$ is a main eigenvector associated with the eigenvalue 0 of $T_{a}^{\sigma^{*}}$. Hence, we have proved that all eigenvalues of $T_{a}^{\sigma^{*}}$ are main.
\end{proof1}
\indent Recall that $K_{t_{1}, \cdots,t_{s}}^{\sigma^{*}}$ denote the signed complete multipartite graph obtained from $K_{t_{1}, \cdots,t_{s}}$ by switching one vertex in each part where $t_{s}\geq 2$ and $s\geq2$. Now we prove all eigenvalues of $K_{t_{1}, \cdots,t_{s}}^{\sigma^{*}}$ are main. The following results are deduced from Theorem 1 (a) and (b) of \cite{ess}. Denote by $\lambda_{\rm min}(\Gamma)$ the smallest eigenvalue of $A(\Gamma)$.
\begin{lemma}(\cite[Theorem 1]{ess})\label{old}
\\
(1) $\lambda_{\rm min}(K_{l_{1}*t_{1},\cdots,l_{s}*t_{s}})\geq-t_{1}$.\\
(2) There is exactly one positive eigenvalue of $K_{l_{1}*t_{1},\cdots,l_{s}*t_{s}}$.\\
(3) There are $s-1$ negative eigenvalues $\lambda_{i}$ where $2\leq i\leq s$, and the following inequality holds,
\begin{align*}
t_{s}<-\lambda_{2}<t_{s-1}<-\lambda_{3}<\cdots<-\lambda_{s}<t_{1}.
\end{align*}
\end{lemma}
\begin{proposition}
All eigenvalues of $K_{t_{1},\cdots,t_{s}}^{\sigma^{*}}$ are main when $t_{s}\geq2$ and $s\geq 2$.
\end{proposition}
\begin{proof1}
The eigenvalues of $K_{t_{1},\cdots,t_{s}}$ are $0^{[n-s]}$ and $\lambda_{1},\lambda_{2},\cdots,\lambda_{s}$ by Lemma \ref{old}, where $\lambda_{i}\notin\{0,-t_{i}\}$. If $\lambda\in\{\lambda_{1},\lambda_{2},\cdots,\lambda_{s}\}$ is an eigenvalue of $K_{t_{1},\cdots,t_{s}}$ and $\bm y$ is an associated eigenvector, then vertices in each part admit the same coordinate in $\bm y$, denoted by $y_{1},y_{2},\cdots,y_{s}$, respectively. Then we have the following equations,
\begin{align}
\left\{
\begin{aligned}
\lambda y_{1}&=t_{2}y_{2}+t_{3}y_{3}+\cdots+t_{s}y_{s},\\
\lambda y_{2}&=t_{1}y_{1}+t_{3}y_{3}+\cdots+t_{s}y_{s},\\
&\cdots\\
\lambda y_{s}&=t_{1}y_{1}+t_{2}y_{2}+\cdots+t_{s-1}y_{s-1}\label{r}.
\end{aligned}
\right.
\end{align}
From Equation (\ref{r}), we deduce $y_{i}=\frac{\lambda-t_{j}}{\lambda-t_{i}}y_{j}$. Hence, $y_{i}\neq 0$ for any $1\leq i\leq s$, otherwise $\bm y=\bm 0$.\\
\indent For the eigenvalue $\lambda\in\{\lambda_{1},\lambda_{2},\cdots,\lambda_{s}\}$, $K_{t_{1},\cdots,t_{s}}^{\sigma^{*}}$ has an eigenvector 
\begin{align*}
\bm\beta(\lambda)=P\bm y=\big[-y_{1},y_{1}\bm{j_{t_{1}-1}},-y_{2},y_{2}\bm{j_{t_{2}-1}},\cdots,-y_{s},y_{s}\bm{j_{t_{s}-1}}\big]^{t}.
\end{align*}
Now we claim that $\bm\beta(\lambda)$ is a main eigenvector of $K_{t_{1},\cdots,t_{s}}^{\sigma^{*}}$ associated with the eigenvalue $\lambda$.\\
\indent Suppose to the contrary we have $\bm\beta(\lambda)^{t}\cdot\bm j=0$ for some $\lambda\in\{\lambda_{1},\lambda_{2},\cdots,\lambda_{s}\}$, which implies
\begin{align}
(t_{1}-2)y_{1}+(t_{2}-2)y_{2}+\cdots+(t_{s}-2)y_{s}=\frac{[\lambda-2(s-1)](y_{1}+y_{2}+\cdots+y_{s})}{s-1}=0\label{rr}.
\end{align}
Furthermore, we have $y_{1}+y_{2}+\cdots+y_{s}=0$ or $\lambda=2(s-1)$. \\
\indent If $y_{1}+y_{2}+\cdots+y_{s}=0$, combining with Equation (\ref{rr}), we have $t_{1}y_{1}+t_{2}y_{2}+\cdots+t_{s}y_{s}=0$ which implies $\bm y=\bm 0$ since $\lambda\neq -t_{i}$. If $\lambda=2(s-1)$, substitute it to Equation (\ref{rr}) and we have
\begin{align}
(t_{1}-2)y_{1}+(t_{2}-2)y_{2}+\cdots+(t_{s}-2)y_{s}=0.\label{rrr}
\end{align}
By Lemma \ref{old} (2), we know $\lambda=2(s-1)$ must be the largest eigenvalue of  $K_{t_{1},\cdots,t_{s}}$. Combining Perron-Frobenius Theorem with $y_{i}\neq 0$, we deduce $y_{i}>0$ for any $1\leq i\leq s$. So when $t_{s}\geq2$, we have $(t_{1}-2)y_{1}+(t_{2}-2)y_{2}+\cdots+(t_{s}-2)y_{s}> 0$, which contradicts (\ref{rrr}). Therefore, $\lambda\neq 2(s-1)$. So $\bm\beta(\lambda)$ is a main eigenvector associated with $\lambda$. \\
\indent It may be checked $\bm{e_{1}}+\bm{e_{2}}$ is a main eigenvector associated with the eigenvalue 0 of $K_{t_{1},\cdots,t_{s}}^{\sigma^{*}}$. Therefore, we have proved that all eigenvalues of $K_{t_{1},\cdots,t_{s}}^{\sigma^{*}}$ are main.
\end{proof1}
\indent For general case, the eigenvalues of the signed complete multipartite graphs obtained by switching one vertex in each part may be not all main. For example, the resulting signed graph obtained from $K_{1*3,2*2}$ has two main eigenvalues and one non-main eigenvalue. Fortunately, we may prove that there always exists a switching $\sigma$ such that all eigenvalues of $K_{l_{1}\ast t_{1},\cdots l_{s}\ast t_{s} }^{\sigma}$ are main, namely, Conjecture 1.1 holds for the complete multipartite graphs. To prove that Conjecture 1.1 holds for the graph $S_{n,k}$ and the complete multipartite graphs, we first present some kinds of technical switching in Section 3.
\section{Auxiliary Results}
\indent Lemma \ref{zero} and Lemma \ref{-1} provide a kind of switching to deal with the eigenvalue 0 and $-1$, respectively. 
\begin{lemma}(\cite[Lemma 3.2]{akb})\label{zero}
Let $G$ be a graph with vertex set $\{v_{1}, \cdots, v_{n}\}$ and $\{v_{1}, \cdots, v_{r}\}$ be a set of vertices with the same neighbors set. If we switch $G$ in vertices $v_{1}, \cdots, v_{t}$, for some $t$, $1\leq t\leq r-1$, then the resulting signed graph has at least $r-1$ independent main eigenvectors for the eigenvalue 0.
\end{lemma} 
\begin{lemma}\label{-1}
Let $G$ be a graph with vertex set $\{v_{1}, \cdots, v_{n}\}$ and $\{v_{1}, \cdots, v_{r}\}$ be a set of vertices with the same closed neighbors set. If we switch $G$ in vertices $v_{1}, \cdots, v_{t}$, for some $t$, $1\leq t\leq r-1$, then the resulting signed graph has at least $r-1$ independent main eigenvectors for the eigenvalue $-1$.
\end{lemma}
\begin{proof1}
It may be checked the following $(r-1)$ vectors
\begin{align*}
\{\bm{e_{i}}+\bm{e_{t+1}}:1\leq i\leq t\}\cup\{\bm{e_{1}}+\bm{e_{t+i}}:2\leq i\leq r-t\}
\end{align*}
are eigenvectors associated with the eigenvalue $-1$ and they are main and independent.
\end{proof1}
\indent Let $\bm \beta=[b_{1},b_{2},\cdots,b_{n}]^{t}$ be a vector, and $\bm{{\beta}_{j_{1},j_{2},\cdots,j_{t}}}$ denote the vector obtained by replacing $b_{j_{i}}$ with $-b_{j_{i}}$ in $\bm\beta$ where $1\leq i\leq t$. If $\bm\beta$ is an eigenvector of $G$ associated with the eigenvalue $\lambda$, then $\bm{\beta_{j_{1},j_{2},\cdots,j_{t}}}$ is an eigenvector of $G^{\sigma}$ associated with the eigenvalue $\lambda$, where $G^{\sigma}$ is obtained from $G$ by switching the vertices corresponding to $j_{1},j_{2},\cdots,j_{t}$. In this paper, we mainly use this argument and the following Lemma \ref{one} and Lemma \ref{same} to modify $\bm\beta$ to be a main eigenvector.
\begin{lemma}(\cite[Lemma 3.3]{akb})\label{one}
Let $\bm{\beta}=[b_{1},\cdots,b_{n}]^{t}$ and $b_{i_{1}},\cdots,b_{i_{k}}$ be non-zero distinct real numbers. Then at most one of the vectors $\{\bm\beta,\bm{\beta_{i_{1}}},\cdots,\bm{\beta_{i_{k}}}\}$ is non-main.
\end{lemma}
\indent If some coordinates of $\bm\beta$ are equal, then we may construct a new vectors sequence, which also contains at most one non-main vector. 
\begin{lemma}\label{same}
Let $\bm\beta=[b_{1},\cdots,b_{n}]^{t}$ and $b_{i_{1}}=b_{i_{2}}=\cdots=b_{i_{k}}$ be non-zero real numbers. Then at most one of the vectors $\{\bm\beta,\bm{\beta_{i_{1}}},\bm{\beta_{i_{1},i_{2}}},\cdots,\bm{\beta_{i_{1},i_{2},\cdots,i_{k}}}\}$ is non-main.
\end{lemma}
\begin{proof1}
Write $\mathcal{S}=\{\bm\beta,\bm{\beta_{i_{1}}},\bm{\beta_{i_{1},i_{2}}},\cdots,\bm{\beta_{i_{1},i_{2},\cdots,i_{k}}}\}$. We will prove the results by showing that any two vectors in $\mathcal{S}$ can not be non-main simultaneously. Suppose to the contrary that $\bm{\beta}$ and $\bm{\beta_{i_{1},i_{2},\cdots,i_{s}}}$ are non-main. That means $\bm{\beta}^{t}\cdot\textbf{\emph{j}}=\bm{\beta_{i_{1},i_{2},\cdots,i_{s}}}^{t}\cdot\textbf{\emph{j}}=0$ for some $1\leq s\leq k$. This implies that $s b_{i_{1}}=0$, which contradicts the assumption $b_{i_{1}}\neq 0$. Suppose $\bm{\beta_{i_{1},i_{2},\cdots,i_{r}}}$ and $\bm{\beta_{i_{1},i_{2},\cdots,i_{s}}}$ are non-main. That means $\bm{\beta_{i_{1},i_{2},\cdots,i_{r}}}^{t}\cdot\textbf{\emph{j}}=\bm{\beta_{i_{1},i_{2},\cdots,i_{s}}}^{t}\cdot\textbf{\emph{j}}=0$ for some $1\leq s<r\leq k$. This implies that $(r-s)b_{i_{1}}=0$, which contradicts the assumption $b_{i_{1}}\neq 0$. 
\end{proof1}
\section{Conjecture \ref{conject} holds for the graph $S_{n,r}$}
\indent Recall graph $S_{n,r}$ is obtained from the complete graph $K_{n-r}$ by attaching $r$ pendent edges at some vertex of $K_{n-r}$. In this section, we prove Conjecture 1.1 holds for $S_{n,r}$. If $n=r+1$ or $n=r+2$, the graph is a star and we know Conjecture 1.1 holds for the complete bipartite graphs from Theorem 4.1 of \cite{akb}. Hence, we suppose $n\geq r+3$. We label the pendent vertices of  $S_{n,r}$ as $v_{1},v_{2},\cdots,v_{r}$, the vertex with degree $n-1$ as $v_{r+1}$, the vertices with degree $n-r-1$ as $v_{r+2},v_{r+3},\cdots,v_{n}$.     
\begin{lemma}\label{component}
(1) The eigenvalues of $S_{n,r}$ $(n\geq r+3)$ are $0^{[r-1]}$, $(-1)^{[n-r-2]}$ and three distinct roots $\lambda_{1},\lambda_{2}$ and $\lambda_{3}$ of the following equation,
\begin{align}
\lambda^{3}-(n-r-2)\lambda^{2}-(n-1)\lambda+r(n-r-2)=0.\label{cha}
\end{align}
(2) If $\lambda\notin\{0,-1\}$ is an eigenvalue of $S_{n,r}$, then $\lambda$ has an eigenvector $\bm x$ with entries given by
\begin{align*}
x_{i}=\left\{
\begin{aligned}
1,\quad\quad \quad&if\ 1\leq i\leq r;\\
\lambda,\quad\quad \quad&if\ i=r+1;\\
\lambda-\frac{r}{\lambda+1},\quad&if\ r+2\leq i\leq n.
\end{aligned}
\right.
\end{align*}
\end{lemma}
\begin{proof1}
Since all the pendent vertices $v_{1},v_{2},\cdots,v_{r}$ have the same neighbors set, by Lemma \ref{zero}, the multiplicity of the eigenvalue 0 is at least $r-1$. The vertices $v_{r+2},v_{r+3},\cdots,v_{n}$ have the same closed neighbors set, by Lemma \ref{-1}, then the multiplicity of the eigenvalue $-1$ is at least $n-r-2$.\\
\indent Let $\bm x$ be an eigenvector associated with the eigenvalue $\lambda\notin\{0,-1\}$. Since $\lambda\neq 0$, all pendent vertices admit the same coordinate in $\bm x$ (denoted by $a$). Similarly, $\lambda\neq -1$ implies that the vertices with degree $n-r-1$ admit the same coordinate in $\bm x$ (denoted by $c$). We assume the vertex with degree $n-1$ admit the coordinate $b$ in $\bm x$. $A(S_{n,r})\bm x=\lambda \bm x$ implies
\begin{align*}
\left\{
\begin{aligned}
&\lambda a=b,\\
&\lambda b=ra+(n-r-1)c,\\
&\lambda c=b+(n-r-2)c.
\end{aligned}
\right.
\end{align*}
Sorting out the above equations, we have 
\begin{align*}
\lambda^{3}-(n-r-2)\lambda^{2}-(n-1)\lambda+r(n-r-2)=0,
\end{align*}
and the associated eigenvector is
\begin{align*}
\bigg[a\bm{j_{r}},\lambda a,(\lambda-\frac{r}{\lambda+1})a\bm{j_{n-r-1}}\bigg]^{t},
\end{align*}
where $\bm{j_{k}}$ is the vector with all-1 entries of order $k$. Taking $a=1$, then we obtain $\bm x=\big[\bm{j_{r}},\lambda,(\lambda-\frac{r}{\lambda+1})\bm{j_{n-r-1}}\big]^{t}$.\\
\indent Let $f(x)=x^{3}-(n-r-2)x^{2}-(n-1)x+r(n-r-2)$ where $n\geq r+3$. Then we have
\begin{align*}
&f(-\infty)<0,\\
&f(0)=r(n-r-2)>0,\\
&f(\sqrt{r})=(-n+1+r)\sqrt{r}<0,\\
&f(+\infty)>0.
\end{align*}
So $f(x)=0$ has three distinct roots which lie in $(-\infty,0)$, $(0,\sqrt{r})$ and $(\sqrt{r},+\infty)$. Hence, the eigenvalues of $S_{n,r}$ $(n\geq r+3)$ are $0^{[r-1]}$, $(-1)^{[n-r-2]}$ and three distinct roots $\lambda_{1},\lambda_{2}$ and $\lambda_{3}$ of Equation (\ref{cha}).
\end{proof1}
\begin{theorem}
There exists a switching $\sigma$ of $S_{n,r}$ such that all eigenvalues of $S_{n,r}^{\sigma}$ are main.
\end{theorem}
\begin{proof1}
From Lemma \ref{component}, $S_{n,r}$ has five distinct eigenvalues, which are $\lambda_{1}$, $\lambda_{2}$, $\lambda_{3}$, $-1$ and 0. Let $P$=diag$(-1,1,\cdots,1,-1)$ and $\bm x$ be defined in Lemma \ref{component}. For $\lambda\in\{\lambda_{1},\lambda_{2},\lambda_{3}\}$, we have $P\bm x$ is an eigenvector of $S_{n,r}^{\sigma_{1}}$, which is obtained from $S_{n,k}$ by switching vertices $v_{1}$ and $v_{n}$. Write $P\bm x=\bm\beta(\lambda)$, namely
\begin{align*}
\bm\beta(\lambda)=\big[-1,\bm{j_{r-1}},\lambda,(\lambda-\frac{r}{\lambda+1})\bm{j_{n-r-2}},-\lambda+\frac{r}{\lambda+1}\big]^{t}.
\end{align*}
\indent If $r\in\{1,2\}$, we will show $\bm\beta(\lambda)$ is main. If $r=1$, suppose to the contrary that $\bm{\beta}(\lambda)$ is a non-main eigenvector corresponding to the eigenvalue $\lambda\in\{\lambda_{1},\lambda_{2},\lambda_{3}\}$, then
\begin{align*}
\bm{\beta}^{t}(\lambda)\cdot\bm{j}=\frac{(n-3)\lambda^{2}+(n-4)\lambda-(n-3)}{\lambda+1}=0.
\end{align*}
Combining with Equation (\ref{cha}), we have $\lambda^{3}-3\lambda=0$. Furthermore, we get $n=3$, $n=2\sqrt{3}$ or $n=-2\sqrt{3}$, which are contradictions. If $r=2$, suppose to the contrary that $\bm\beta^{t}(\lambda)\cdot\bm j=0$, namely 
\begin{align*}
\bm{\beta}^{t}(\lambda)\cdot\bm{j}=\frac{(n-4)\lambda^{2}+(n-4)\lambda-2(n-5)}{\lambda+1}=0.
\end{align*}
Combining with Equation (\ref{cha}), we have $\lambda^{3}-3\lambda+2=0$, which implies $\lambda=1$ or $\lambda=-2$, while it contradicts Equation (\ref{cha}).\\
\indent  If $r\geq3$ and $n=r+3$, we claim $\bm{\beta}(\lambda)$ is main, otherwise
\begin{align*}
\bm{\beta}(\lambda)^{t}\cdot \textbf{\emph{j}}=r-2+\lambda=0.
\end{align*}
Substituting $\lambda=2-r$ into Equation (\ref{cha}), we obtain
\begin{align*}
r(r^{2}-6r+7)=0,
\end{align*}
which contradicts the fact that $r$ is an integer.  \\
\indent It may be checked that $\bm{e_{1}}+\bm{e_{2}}$ is a main eigenvector associated with the eigenvalue 0 of $S_{n,k}^{\sigma_{1}}$ and $\bm{e_{n}}+\bm{e_{n-1}}$ is a main eigenvector associated with the eigenvalue $-1$ of $S_{n,k}^{\sigma_{1}}$. Therefore, we have proved that all eigenvalues of $S_{n,r}^{\sigma_{1}}$ are main when $r$=1, 2, or $r\geq 3$ and $n=r+3$.\\
\indent Now, we suppose $r\geq 3$ and $n\geq r+4$.  If $\lambda-\frac{r}{\lambda+1}=0$, combining with Equation (\ref{cha}), we get $r=0$. If $1=\lambda-\frac{r}{\lambda+1}$, combining with Equation (\ref{cha}), we get $\lambda=-1$. If $\lambda=1$, according to Equation (\ref{cha}), we have 
\begin{align*}
\left\{
\begin{aligned}
&(1-r)b=(n-r-1)c,\\
&-b=(n-r-3)c,
\end{aligned}
\right.
\end{align*}
then
\begin{align*}
n=r+3+\frac{2}{r-2}.
\end{align*}
We can easily see that $n$ is not an integer if $r\geq 4$. When $r=3$, we have $n=7$, while $\lambda=1$ is not an eigenvalue of $S_{7,3}$. So when $\lambda\in\{\lambda_{1},\lambda_{2},\lambda_{3}\}$, all elements in $\{1,\lambda,\lambda-\frac{r}{\lambda+1}\}$ are non-zero distinct real numbers. Set
\begin{align*}
\mathcal{S}=\big\{\bm{\beta}(\lambda),\bm{\beta_{2}}(\lambda),\bm{\beta_{r+1}}(\lambda),\bm{\beta_{n-1}}(\lambda)\big\}.
\end{align*}
By Lemma \ref{one}, for a fixed $\lambda_{i}$, at most one vector in $\mathcal{S}$ is non-main. Noting $|\mathcal{S}|=4$, hence there exists a vector $\bm\gamma(\lambda)\in\mathcal{S}$, which is main for every $\lambda\in\{\lambda_{1},\lambda_{2},\lambda_{3}\}$. We suppose $\sigma_{2}$ is the switching corresponding to $\bm\gamma(\lambda)$. In $\sigma_{2}$, we suppose that $p$ $(1\leq p\leq 2)$ pendent vertices will be switched and $q$ $(1\leq q\leq 2)$ vertices with degree $n-r-1$ will be switched. It may be checked that $\bm{e_{1}}+\bm{e_{p+1}}$ is a main eigenvector of the eigenvalue 0 of $S_{n,k}^{\sigma_{2}}$ and $\bm{e_{n}}+\bm{e_{n-q}}$ is a main eigenvector of the eigenvalue $-1$ of $S_{n,k}^{\sigma_{2}}$. Therefore, we have proved that all eigenvalues of $S_{n,r}^{\sigma_{2}}$ are main.
\end{proof1}
\section{Conjecture \ref{conject} holds for complete multipartite graphs}
\indent We will prove Conjecture 1.1 holds for complete multipartite graphs in this section. It was pointed out in \cite{akb} that Conjecture 1.1 holds for all graphs up to 9 vertices by a computer search. Hence, we only need to consider complete multipartite graphs with not less than 9 vertices. Set $V(K_{l_{1}*t_{1},\cdots,l_{s}*t_{s}})=\mathop{\cup}\limits_{i=1}^{s}U_{i}$ with $|U_{i}|=m_{i}=l_{i}t_{i}$ and $U_{i}=\mathop{\cup}\limits_{j=1}^{l_{i}}U_{i,j}$ with $|U_{i,j}|=t_{i}$, $1\leq j\leq l_{i}$, namely, $U_{i,1},U_{i,2},\cdots,U_{i,l_{i}}$ are $l_{i}$ parts of $K_{l_{1}*t_{1},\cdots,l_{s}*t_{s}}$. We will label the vertices from $U_{1,1}$ to $U_{s,l_{s}}$ as $v_{1},v_{2},\cdots,v_{n}$. Write $f_{1}=0$, and $f_{i}=\sum\limits_{j=1}^{i-1}m_{j}$. 
\begin{lemma}\label{-n}
(1) $-t_{i}$ is an eigenvalue of $K_{l_{1}*t_{1},\cdots,l_{s}*t_{s}}$ with multiplicity at least $l_{i}-1$.\\
(2) If we switch $p$ vertices in some $U_{i,j}$ where $1\leq p\leq t_{i}$, $1\leq j\leq l_{i}$ and switch $q$ vertices in some $U_{i,k}$ where $0\leq q<p$, $1\leq k\neq j\leq l_{i}$, then the eigenvalue $-t_{i}$ of the resulting signed graph is main.
\end{lemma}
\begin{proof1}
(1) It is not difficult to see the following $(l_{i}-1)$ linearly independent vectors
\begin{align}
\bm{z_{j}}=\sum\limits_{v_{k}\in U_{i,1}}\bm{e_{k}}-\sum\limits_{v_{k}\in U_{i,j}}\bm{e_{k}},\quad 2\leq j\leq l_{i}\label{-ni},
\end{align}
are eigenvectors associated with the eigenvalue $-t_{i}$ of $K_{l_{1}*t_{1},\cdots,l_{s}*t_{s}}$.\\
(2) Without loss of generality, we switch $v_{f_{i}+1},\cdots,v_{f_{i}+p}$ in $U_{i,1}$ where $1\leq p\leq t_{i}$ and switch $v_{f_{i}+t_{i}+1}, v_{f_{i}+t_{i}+2},\cdots,v_{f_{i}+t_{i}+q}$ in $U_{i,2}$ where $0\leq q<p$. The resulting signed graph is denoted by $K_{l_{1}*t_{1},\cdots,l_{s}*t_{s}}^{\sigma}$. By (\ref{-ni}),
\begin{align*}
\bm{z_{2}^{'}}=-\sum\limits_{j=1}^{p}\bm{e_{f_{i}+j}}+\sum\limits_{j=p+1}^{t_{i}}\bm{e_{f_{i}+j}}+\sum\limits_{j=1}^{q}\bm{e_{f_{i}+t_{i}+j}}-\sum\limits_{j=q+1}^{t_{i}}\bm{e_{f_{i}+t_{i}+j}}
\end{align*}
is an eigenvector associated with the eigenvalue $-t_{i}$ of $K_{l_{1}*t_{1},\cdots,l_{s}*t_{s}}^{\sigma}$. Furthermore $\bm{z_{2}^{'}}$ is main, since $(\bm{z_{2}^{'}})^{t}\cdot\bm j=2(q-p)\neq 0$.  
\end{proof1}
\indent The multiplicity of the eigenvalue 0 of $K_{l_{1}*t_{1},\cdots,l_{s}*t_{s}}$ is at least $b=\sum\limits_{i=1}^{s}(m_{i}-l_{i})$ (see \cite[Lemma 2]{ess}). By Lemma \ref{-n}, the multiplicity of the eigenvalue $-t_{i}$ of $K_{l_{1}*t_{1},\cdots,l_{s}*t_{s}}$ is at least $l_{i}-1$. Moreover, by Lemma \ref{old} (2) (3), we deduce that there are exactly $s$ distinct eigenvalues of $K_{l_{1}*t_{1},\cdots,l_{s}*t_{s}}$ different from $-t_{i}$ and 0, say $\lambda_{1},\lambda_{2},\cdots,\lambda_{s}$. Hence, the eigenvalues of $K_{l_{1}*t_{1},\cdots,l_{s}*t_{s}}$ are $0^{[b]}$, $(-t_{i})^{[l_{i}-1]}$, $\lambda_{1}$, $\lambda_{2}$, $\cdots$, $\lambda_{s}$. \\
\indent If $\lambda\in\{\lambda_{1},\lambda_{2},\cdots,\lambda_{s}\}$ and $\bm z$ is an associated eigenvector, then the vertices in $U_{i}$ admit the same coordinate in $\bm z$, denoted by $z_{i}$, $1\leq i\leq s$. Then
\begin{align}
\left\{
 \begin{aligned}
(\lambda+t_{1}) z_{1}&=m_{1}z_{1}+m_{2}z_{2}+\cdots+m_{s}z_{s},\\
(\lambda+t_{2}) z_{2}&=m_{1}z_{1}+m_{2}z_{2}+\cdots+m_{s}z_{s},\\
&\cdots\\
(\lambda+t_{s}) z_{s}&=m_{1}z_{1}+m_{2}z_{2}+\cdots+m_{s}z_{s}\label{s}.
\end{aligned}
\right.
\end{align}
From Equation (\ref{s}), we deduce $z_{i}=\frac{\lambda+t_{j}}{\lambda+t_{i}}z_{j}$. Hence $z_{i}\neq 0$ for any $1\leq i\leq s$, otherwise $\bm z=\bm 0$. Moreover, the assumption $t_{i}\neq t_{j}$ implies $z_{i}\neq z_{j}$.\\
\indent Three rules should be obeyed during the later switching operations.\\
\indent (R1) During the whole operations, each vertex is allowed to be switched at most once.\\
\indent (R2) When $t_{i}=3$ and $m_{i}\geq6$, if we need to switch three vertices in $U_{i}$, then the vertices $v_{f_{i}+1}$, $v_{f_{i}+2}$ and $v_{f_{i}+4}$ (not $v_{f_{i}+1}$, $v_{f_{i}+2}$ and $v_{f_{i}+3}$) should be picked. \\
\indent (R3) After rules R2 and R3 are met, we always choose the vertex with the smallest labelling if we need to switch more vertices. \\
\indent We use $\bm{y_{i^{[j]}}}$ to denote the vector obtained from $\bm y$ by switching $j$ vertices in $U_{i}$. 
\begin{lemma}\label{1}
If $t_{s}\geq 2$, $K_{l_{1}*t_{1},\cdots,l_{s}*t_{s}}$ has a switching $\sigma$ such that all eigenvalues of $K_{l_{1}*t_{1},\cdots,l_{s}*t_{s}}^{\sigma}$ are main.
\end{lemma}
\begin{proof1}
If $s=1$ and $l_{s}=1$, $K_{l_{1}*t_{1},\cdots,l_{s}*t_{s}}$ is an empty graph, and the result is trivial. \\
\indent If $s=1$ and $l_{s}\geq 2$, we know that the eigenvalues of $K_{l_{1}*t_{1}}$ are $(-t_{1})^{[l_{1}-1]},0^{[n-l_{1}]},(l_{1}-1)t_{1}$. The signed graph after switching $v_{1}$ of $K_{l_{1}*t_{1}}$ is denoted by $K_{l_{1}*t_{1}}^{\sigma_{1}}$, and then the eigenvalue $(l_{1}-1)t_{1}$ of $K_{l_{1}*t_{1}}^{\sigma_{1}}$ is main with eigenvector $[-1,1,\cdots,1]^{t}$. By using the arguments of Lemma \ref{zero} and Lemma \ref{-n}, we deduce $\bm{e_{1}}+\bm{e_{2}}$ is a main eigenvector associated with the eigenvalue 0 of $K_{l_{1}*t_{1}}^{\sigma_{1}}$ and 
\begin{align*}
-\bm{e_{1}}+\sum\limits_{j=2}^{t_{1}}\bm{e_{j}}-\sum\limits_{j=1}^{t_{1}}\bm{e_{t_{1}+j}}
\end{align*}
is a main eigenvector associated with the eigenvalue $-t_{1}$ of $K_{l_{1}*t_{1}}^{\sigma_{1}}$.\\
\indent If $s\geq 2$, then we have $t_{s-1}\geq 3$. If $l_{i}\geq 2$ ($1\leq i\leq s-1$), then we will switch the vertex with the smallest labelling in $U_{i}$. The new eigenvector associated with the eigenvalue $\lambda\in\{\lambda_{1},\lambda_{2},\cdots,\lambda_{s}\}$ after this switching and switching $v_{f_{s}+1}$ is denoted by $\bm\eta(\lambda)$. Set
\begin{align*}
\mathcal{F}=\big\{\bm{\eta}(\lambda),\bm{\eta_{1^{[1]}}}(\lambda),\bm{\eta_{1^{[2]}}}(\lambda),\bm{\eta_{2^{[2]}}}(\lambda),\cdots,\bm{\eta_{(s-1)^{[2]}}}(\lambda)\big\}.
\end{align*}
\indent For a fixed $\lambda_{i}$, we will prove that at most one vector in $\mathcal{F}$ is non-main. Assume $\bm{\eta}(\lambda)$ and $\bm{\eta_{1^{[1]}}}(\lambda)$ are non-main. That means $\bm{\eta}(\lambda)^{t}\cdot\textbf{\emph j}=\bm{\eta_{1^{[1]}}}(\lambda)^{t}\cdot\textbf{\emph j}=0$, which implies a contradiction $z_{1}=0$. Similarly, $\bm{\eta}(\lambda)$ and $\bm{\eta_{i^{[2]}}}(\lambda)$ can not be non-main simultaneously. Assume $\bm{\eta_{i^{[2]}}}(\lambda)$ and $\bm{\eta_{j^{[2]}}}(\lambda)$ are non-main where $i\neq j$, which implies a contradiction $z_{i}=z_{j}$. Suppose $\bm{\eta_{1^{[1]}}}(\lambda)$ and $\bm{\eta_{i^{[2]}}}(\lambda)$ are non-main, which implies $z_{1}=2z_{i}$. Combining with Equation (\ref{s}), we have $\frac{1}{\lambda+t_{1}}=\frac{2}{\lambda+t_{i}}$, thus $\lambda=-2t_{1}+t_{i}<-t_{1}$. While from Lemma \ref{old} (1) we know $\lambda\geq\lambda_{\rm min}(K_{l_{1}*t_{1},\cdots,l_{s}*t_{s}})>-t_{1}$.\\
\indent Since for a fixed $\lambda_{i}$ at most one vector in $\mathcal{F}$ is non-main and $|\mathcal{F}|=s+1$, there exists a main eigenvector $\bm\gamma(\lambda)\in\mathcal{F}$ for every $\lambda\in\{\lambda_{1},\lambda_{2},\cdots,\lambda_{s}\}$. The switching corresponding to $\bm\gamma(\lambda)$ is denoted by $\sigma_{2}$. Then $\bm\gamma(\lambda)$ is a main eigenvector of $K_{l_{1}*t_{1},\cdots,l_{s}*t_{s}}^{\sigma_{2}}$ associated with $\lambda$. Now, we need to show the other eigenvalues of $K_{l_{1}*t_{1},\cdots,l_{s}*t_{s}}^{\sigma_{2}}$ are also main. Since in $\sigma_{2}$ only $v_{f_{s}+1}$ has been switched in $U_{s}$ and $t_{s}\geq 2$, by Lemma \ref{zero}, $\bm{e_{f_{s}+1}}+\bm{e_{f_{s}+2}}$ is a main eigenvector associated with the eigenvalue 0 of $K_{l_{1}*t_{1},\cdots,l_{s}*t_{s}}^{\sigma_{2}}$. If $-t_{i}$ is an eigenvalue of $K_{l_{1}*t_{1},\cdots,l_{s}*t_{s}}^{\sigma_{2}}$, then $l_{i}\geq 2$. Hence there are $w\in\{1,2,3\}$ vertices switched in $U_{i}$ in $\sigma_{2}$. Accoording to the switching operation rules, if $w\in\{1,2\}$, there are $w$ vertices switched in $U_{i,1}$ and no vertices switched in $U_{i,2}$. If $w=3$ and $t_{i}=2$, there are two vertices switched in $U_{i,1}$ and one vertices switched in $U_{i,2}$. If $w=3$ and $t_{i}=3$, combining with the Rule 2, there are two vertices switched in $U_{i,1}$ and one vertices switched in $U_{i,2}$. If $w=3$ and $t_{i}\geq 4$, there are three vertices switched in $U_{i,1}$ and no vertices switched in $U_{i,2}$. In general, there are $1\leq p\leq t_{i}$ vertices switched in $U_{i,1}$ and $0\leq q<p$ vertices switched in $U_{i,2}$. Hence, by Lemma \ref{-n}, the following vector
\begin{align*}
-\sum\limits_{j=1}^{p}\bm{e_{f_{i}+j}}+\sum\limits_{j=p+1}^{t_{i}}\bm{e_{f_{i}+j}}+\sum\limits_{j=1}^{q}\bm{e_{f_{i}+t_{i}+j}}-\sum\limits_{j=q+1}^{t_{i}}\bm{e_{f_{i}+t_{i}+j}}
\end{align*}
is a main eigenvector of the eigenvalue $-t_{i}$ of $K_{l_{1}*t_{1},\cdots,l_{s}*t_{s}}^{\sigma_{2}}$. We have proved all eigenvalues of $K_{l_{1}*t_{1},\cdots,l_{s}*t_{s}}^{\sigma_{2}}$ are main.
\end{proof1}
\begin{lemma}\label{2}
If $t_{s}=1$, $s\geq2$ and $l_{i}=1$, $1\leq i\leq s-1$, $K_{l_{1}*t_{1},\cdots,l_{s}*t_{s}}$ has a switching $\sigma$ such that all eigenvalues of $K_{l_{1}*t_{1},\cdots,l_{s}*t_{s}}^{\sigma}$ are main.
\end{lemma}
\begin{proof1}
If $t_{1}=2$, then $m_{1}=2$, $s=2$, and $t_{2}=1$. If $l_{2}\leq3$, then the number of vertices of $K_{l_{1}*t_{1},\cdots,l_{s}*t_{s}}$ is at most 5. If $l_{2}\geq 4$, set $\bm{\xi}(\lambda)=\big[-z_{1},z_{1}\bm{j_{m_{1}-1}},-z_{2},z_{2}\bm {j_{m_{2}-1}}\big]^{t}$ and 
\begin{align*}
\mathcal{F}_{1}=\big\{\bm{\xi}(\lambda),\bm{\xi_{m_{1}+2}}(\lambda),\bm{\xi_{m_{1}+2,m_{1}+3}}(\lambda)\big\}.
\end{align*}
Since $v_{m_{1}+2}$ and $v_{m_{1}+3}$ admit the same coordinate in $\bm{\xi}(\lambda)$, by Lemma \ref{same}, for a fixed $\lambda_{i}$, at most one vector in $\mathcal{F}_{1}$ is non-main. Noting $|\mathcal{F}_{1}|=3$, there exists a main eigenvector $\bm{\gamma_{1}}(\lambda)\in\mathcal{F}_{1}$ for every $\lambda\in\{\lambda_{1},\lambda_{2}\}$. The switching corresponding to $\bm{\gamma_{1}}(\lambda)$ is denoted by $\sigma_{1}$ which means switching $v_{1}$ in $U_{1}$ and $q\in\{1,2,3\}$ vertices in $U_{2}$. By using the arguments of Lemma \ref{zero} and Lemma \ref{-1}, we know $\bm{e_{1}}+\bm{e_{2}}$ is a main eigenvectors associated with the eigenvalue $0$ of $K_{l_{1}*t_{1},\cdots,l_{s}*t_{s}}^{\sigma_{1}}$ and $\bm {e_{m_{1}+1}}+\bm{e_{m_{1}+q+1}}$is a main eigenvectors associated with the eigenvalue $-1$ of $K^{\sigma_{1}}_{l_{1}*t_{1},\cdots,l_{s}*t_{s}}$. Therefore, all eigenvalues of $K_{l_{1}*t_{1},\cdots,l_{s}*t_{s}}^{\sigma_{1}}$ are main.\\
\indent If $t_{1}=3$ and $l_{s}\leq 2$, the number of vertices of $K_{l_{1}*t_{1},\cdots,l_{s}*t_{s}}$ is at most 7. If $l_{s}\geq3$, then $m_{s}\geq 3$ and let 
\begin{align*}
\bm\zeta(\lambda)=\big[-z_{1},z_{1}\bm{j_{m_{1}-1}},z_{2}\bm{j_{m_{2}}},\cdots,z_{s-1}\bm{j_{m_{s-1}}},-z_{s},z_{s}\bm{j_{m_{s}-1}}\big]^{t}
\end{align*} 
be the eigenvector associated with $\lambda\in\{\lambda_{1},\lambda_{2},\cdots,\lambda_{s}\}$ after switching vertices $v_{1}$ and $v_{f_{s}+1}$. We consider
\begin{align*}
\mathcal{F}_{2}=\big\{\bm{\zeta}(\lambda),\bm{\zeta_{1^{[1]}}}(\lambda),\cdots,\bm{\zeta_{s^{[1]}}}(\lambda)\big\}.
\end{align*}
Since $z_{i}\neq z_{j}$, $i\neq j$, by Lemma \ref{one}, for a fixed $\lambda_{i}$, at most one vector in $\mathcal{F}_{2}$ is non-main. Noting $|\mathcal{F}_{2}|=s+1$, there exists a main eigenvector $\bm{\gamma_{2}}(\lambda)\in \mathcal{F}_{2}$ for every $\lambda\in\{\lambda_{1},\lambda_{2},\cdots,\lambda_{s}\}$. The switching corresponding to $\bm{\gamma_{2}}(\lambda)$ is denoted by $\sigma_{2}$. Since $w\in\{1,2\}$ vertices have been switched in $U_{s}$ in $\sigma_{2}$ and $m_{s}\geq 3$, by Lemma \ref{-1}, $\bm{e_{f_{s}+1}}+\bm{e_{f_{s}+w+1}}$ is a main eigenvector associated with the eigenvalue $-1$ of $K_{l_{1}*t_{1},\cdots,l_{s}*t_{s}}^{\sigma_{2}}$. Similarly, since $w\in\{1,2\}$ vertices have been switched in $U_{1}$ in $\sigma_{2}$ and $t_{1}=3$, by Lemma \ref{zero}, $\bm{e_{1}}+\bm{e_{w+1}}$ is a main eigenvector associated with the eigenvalue 0 of $K_{l_{1}*t_{1},\cdots,l_{s}*t_{s}}^{\sigma_{2}}$.\\
\indent If $t_{1}\geq4$, set
\begin{align*}
\mathcal{F}_{3}=\big\{\bm{\zeta}(\lambda),\bm{\zeta_{1^{[1]}}}(\lambda),\bm{\zeta_{1^{[2]}}}(\lambda),\bm{\zeta_{2^{[2]}}}(\lambda),\cdots,\bm{\zeta_{(s-1)^{[2]}}}(\lambda)\big\}.
\end{align*}
By using the similar arguments as the vectors in $\mathcal{F}$ in Lemma \ref{1}, we may prove that for a fixed $\lambda_{i}$, at most one vector in $\mathcal{F}_{3}$ is non-main. Noting $|\mathcal{F}_{3}|=s+1$, then there exists a main eigenvector $\bm{\gamma_{3}}(\lambda)\in\mathcal{F}_{3}$ for every $\lambda\in\{\lambda_{1},\lambda_{2},\cdots,\lambda_{s}\}$. The switching corresponding to $\bm{\gamma_{3}}(\lambda)$ is denoted by $\sigma_{3}$. Since only $v_{f_{s}+1}$ has been switched in $U_{s}$ in $\sigma_{3}$, by Lemma \ref{-1}, $\bm{e_{f_{s}+1}}+\bm{e_{f_{s}+2}}$ is a main eigenvector of the eigenvalue $-1$ of $K_{l_{1}*t_{1},\cdots,l_{s}*t_{s}}^{\sigma_{3}}$ when $l_{s}\geq 2$. Since $w\in\{1,2,3\}$ vertices have been switched in $U_{1}$ in $\sigma_{3}$ and $t_{1}\geq4$, $\bm{e_{1}}+\bm{e_{w+1}}$ is a main eigenvector of the eigenvalue 0 of $K_{l_{1}*t_{1},\cdots,l_{s}*t_{s}}^{\sigma_{3}}$ by Lemma \ref{zero}. Therefore, all eigenvalues of $K_{l_{1}*t_{1},\cdots,l_{s}*t_{s}}^{\sigma_{3}}$ are main.
\end{proof1}
\begin{lemma}\label{3}
If $t_{s}=1$, $s\geq2$ and there exists some $l_{i}\geq 2$, $1\leq i\leq s-1$, $K_{l_{1}*t_{1},\cdots,l_{s}*t_{s}}$ has a switching $\sigma$ such that all eigenvalues of $K_{l_{1}*t_{1},\cdots,l_{s}*t_{s}}^{\sigma}$ are main.
\end{lemma}
\begin{proof1}
If $t_{1}=2$, then $s=2$. From the assumption that there exists some $l_{i}\geq 2$ $(1\leq i\leq s-1)$ we deduce $l_{1}\geq 2$. Furthermore, $m_{1}=l_{1}t_{1}\geq 2t_{1}=4$. If $l_{2}\geq 2$, all distinct eigenvalues of $K_{l_{1}*t_{1},l_{2}*t_{2}}$ are $-2$, $-1$, 0, $\lambda_{1}$, $\lambda_{2}$ and if $l_{2}=1$, all distinct eigenvalues are $-2$, 0, $\lambda_{1}$, $\lambda_{2}$. \\
\indent If $1\leq m_{2}\leq 3$ and $m_{1}=4$, the number of vertices of $K_{l_{1}*t_{1},l_{2}*t_{2}}$ is at most 7. If $1\leq m_{2}\leq 3$ and $m_{1}\geq 6$, set $\bm{\xi}(\lambda)=\big[-z_{1},z_{1}\bm{j_{m_{1}-1}},-z_{2},z_{2}\bm{j_{m_{2}-1}}\big]^{t}$ and 
\begin{align*}
\mathcal{F}_{1}=\big\{\bm{\xi}(\lambda),\bm{\xi_{2,3}}(\lambda),\bm{\xi_{2,3,4,5}}(\lambda)\big\}.
\end{align*}
Since $v_{2},v_{3},v_{4}$ and $v_{5}$ admit the same coordinate in $\bm{\xi}(\lambda)$ and $\mathcal{F}_{1}\subseteq\{\bm{\xi}(\lambda),\bm{\xi_{2}}(\lambda),\bm{\xi_{2,3}}(\lambda),\bm{\xi_{2,3,4}}(\lambda),$\\
$\bm{\xi_{2,3,4,5}}(\lambda)\}$, by Lemma \ref{same}, for a fixed $\lambda_{i}$, at most one vector in $\mathcal{F}_{1}$ is non-main. Noting $|\mathcal{F}_{1}|=3$, there exists a main eigenvector $\bm{\gamma_{1}}(\lambda)\in\mathcal{F}_{1}$ for every $\lambda\in\{\lambda_{1},\lambda_{2}\}$. \\
\indent The switching corresponding to $\bm{\gamma_{1}}(\lambda)$ is denoted by $\sigma_{1}$. In $\sigma_{1}$, $v_{f_{2}+1}$ has been switched in $U_{2}$ and $q\in\{1,3,5\}$ vertices have been switched in $U_{1}$. In fact, $\bm e_{q}+\bm e_{q+1}$ is a main eigenvector associated with the eigenvalue 0 of $K_{l_{1}*t_{1},l_{2}*t_{2}}^{\sigma_{1}}$. If $l_{2}\geq 2$, $\bm{e_{m_{1}+1}}+\bm{e_{m_{1}+2}}$ is a main eigenvector associated with the eigenvalue $-1$ of $K_{l_{1}*t_{1},l_{2}*t_{2}}^{\sigma_{1}}$. If $q\in\{1,3\}$, then $-\bm{e_{q}}+\bm{e_{q+1}}-\bm{e_{q+2}}-\bm{e_{q+3}}$ is a main eigenvector associated with the eigenvalue $-2$ of $K_{l_{1}*t_{1},l_{2}*t_{2}}^{\sigma_{1}}$. If $q=5$, then $-\bm{e_{1}}-\bm{e_{2}}+\bm{e_{5}}-\bm{e_{6}}$ is a main eigenvector associated with the eigenvalue $-2$ of $K_{l_{1}*t_{1},l_{2}*t_{2}}^{\sigma_{1}}$. Therefore, all eigenvalues of $K_{l_{1}*t_{1},l_{2}*t_{2}}^{\sigma_{1}}$ are main.\\
\indent Suppose $m_{2}\geq 4$ and set
\begin{align*}
\mathcal{F}_{2}=\big\{\bm{\xi}(\lambda),\bm{\xi_{m_{1}+2}}(\lambda),\bm{\xi_{m_{1}+2,m_{1}+3}}(\lambda)\big\}.
\end{align*}
Since $v_{m_{1}+2}$ and $v_{m_{1}+3}$ admit the same coordinate in $\bm{\xi}(\lambda)$, by Lemma \ref{same}, for a fixed $\lambda_{i}$, at most one vector in $\mathcal{F}_{2}$ is non-main. Noting $|\mathcal{F}_{2}|=3$, there exists a main eigenvector $\bm{\gamma_{2}}(\lambda)\in\mathcal{F}_{2}$ for every $\lambda\in\{\lambda_{1},\lambda_{2}\}$. The switching corresponding to $\bm{\gamma_{2}}(\lambda)$ is denoted by $\sigma_{2}$ which means switching $v_{1}$ in $U_{1}$ and $q\in\{1,2,3\}$ vertices in $U_{2}$. By using the arguments of Lemma \ref{zero}, Lemma \ref{-1} and Lemma \ref{-n}, we know $\bm{e_{1}}+\bm{e_{2}}$ is a main eigenvectors associated with the eigenvalue $0$ of $K_{l_{1}*t_{1},l_{2}*t_{2}}^{\sigma_{2}}$ and $\bm{e_{m_{1}+1}}+\bm{e_{m_{1}+q+1}}$ is a main eigenvectors associated with the eigenvalue $-1$ of $K_{l_{1}*t_{1},l_{2}*t_{2}}^{\sigma_{2}}$ and $-\bm{e_{1}}+\bm{e_{2}}-\bm{e_{3}}-\bm{e_{4}}$ is a main eigenvectors associated with the eigenvalue $-2$ of $K_{l_{1}*t_{1},l_{2}*t_{2}}^{\sigma_{2}}$. Therefore, all eigenvalues of $K_{l_{1}*t_{1},l_{2}*t_{2}}^{\sigma_{2}}$ are main.\\
\indent Now we need to consider $t_{1}\geq 3$. If $l_{i}\geq 2$ ($1\leq i\leq s-1$), then we will switch the vertex with the smallest labelling in $U_{i}$. The new eigenvector associated with the eigenvalue $\lambda\in\{\lambda_{1},\lambda_{2},\cdots,\lambda_{s}\}$ after this switching and switching $v_{f_{s}+1}$ is denoted by $\bm\eta(\lambda)$. Noting $t_{i}\geq 2$ always holds for $1\leq i\leq s-1$, we set 
\begin{align*}
\mathcal{F}_{3}=\big\{\bm{\eta}(\lambda),\bm{\eta_{1^{[1]}}}(\lambda),\bm{\eta_{1^{[2]}}}(\lambda),\bm{\eta_{2^{[2]}}}(\lambda),\cdots,\bm{\eta_{(s-1)^{[2]}}}(\lambda)\big\}.
\end{align*}
By using the similar arguments as the vectors in $\mathcal{F}$ in Lemma \ref{1}, we may prove that for a fixed $\lambda_{i}$, at most one vector in $\mathcal{F}_{3}$ is non-main. Noting $|\mathcal{F}_{3}|=s+1$, there exists a main eigenvector $\bm{\gamma_{3}}(\lambda)\in\mathcal{F}_{3}$ for every $\lambda\in\{\lambda_{1},\lambda_{2},\cdots,\lambda_{s}\}$. The switching corresponding to $\bm{\gamma_{3}}(\lambda)$ is denoted by $\sigma_{3}$. Next, we will show the eigenvalues $-t_{i}$, 0 are main eigenvalues of $K_{l_{1}*t_{1},\cdots,l_{s}*t_{s}}^{\sigma_{3}}$. If $-t_{i}$ $(1\leq i\leq s-1)$ is an eigenvalue of $K_{l_{1}*t_{1},\cdots,l_{s}*t_{s}}^{\sigma_{3}}$, then $l_{i}\geq 2$. Hence there are $w\in\{1,2,3\}$ vertices switched in $U_{i}$ in $\sigma_{3}$. Combining $t_{i}\geq 2$ with the switching operation rules, actually there are always $1\leq p\leq t_{i}$ vertices switched in $U_{i,1}$ and $0\leq q<p$ vertices switched in $U_{i,2}$. By Lemma \ref{-n}, the following vector
\begin{align*}
-\sum\limits_{j=1}^{p}\bm{e_{f_{i}+j}}+\sum\limits_{j=p+1}^{t_{i}}\bm{e_{f_{i}+j}}+\sum\limits_{j=1}^{q}\bm{e_{f_{i}+t_{i}+j}}-\sum\limits_{j=q+1}^{t_{i}}\bm{e_{f_{i}+t_{i}+j}}
\end{align*}
is a main eigenvector of the eigenvalue $-t_{i}$ of $K_{l_{1}*t_{1},\cdots,l_{s}*t_{s}}^{\sigma_{3}}$. Since only $v_{f_{s}+1}$ has been switched in $U_{s}$ in $\sigma_{3}$, by Lemma \ref{-1}, $\bm{e_{f_{s}+1}}+\bm{e_{f_{s}+2}}$ is a main eigenvector of the eigenvalue $-1$ of $K_{l_{1}*t_{1},\cdots,l_{s}*t_{s}}^{\sigma_{3}}$ when $l_{s}\geq 2$. \\
\indent If $l_{1}\geq 2$, then $w\in\{1,2,3\}$ vertices have been switched in $U_{1}$ in $\sigma_{3}$. If $l_{1}=1$, then $w\in\{1,2\}$ vertices have been switched in $U_{1}$ in $\sigma_{3}$. Combining $t_{1}\geq 3$ with the rule R3, we find that in either case, $t$ vertices have been switched in $U_{1,1}$ in $\sigma_{3}$ where $1\leq t\leq t_{1}-1$. Hence, by Lemma \ref{zero}, $\bm{e_{1}}+\bm{e_{t+1}}$ is a main eigenvector of the eigenvalue 0 of $K_{l_{1}*t_{1},\cdots,l_{s}*t_{s}}^{\sigma_{3}}$. In this way, we have proved that all eigenvalues of $K_{l_{1}*t_{1},\cdots,l_{s}*t_{s}}^{\sigma_{3}}$ are main.
\end{proof1}
\begin{theorem}
The complete multipartite graph $K_{l_{1}*t_{1},\cdots,l_{s}*t_{s}}$ has a switching $\sigma$ such that all eigenvalues of $K^{\sigma}_{l_{1}*t_{1},\cdots,l_{s}*t_{s}}$ are main.
\end{theorem}
\begin{proof1}
If $t_{s}\geq 2$, the result follows from Lemma \ref{1}.\\
\indent If $t_{s}=1$ and $s=1$, then $K_{l_{1}*t_{1},\cdots,l_{s}*t_{s}}$ is the complete graph with at least 9 vertices, and all eigenvalues are main after switching one vertex. \\
\indent If $t_{s}=1$, $s\geq2$ and $l_{i}=1$, $1\leq i\leq s-1$, the result follows from Lemma \ref{2}.\\
\indent If $t_{s}=1$, $s\geq2$ and there exists some $l_{i}\geq 2$ $(1\leq i\leq s-1)$, the result follows from Lemma \ref{3}.
\end{proof1}
\section*{Acknowledgement}
The authors would like to express their sincere gratitude to Prof. Wei Wang for his careful reading and valuable suggestions.

\end{document}